\documentclass{amsart}

\usepackage[foot]{amsaddr}

\usepackage[T1]{fontenc}
\usepackage[utf8]{inputenc}
\usepackage[USenglish]{babel}
\usepackage{textcase}

\usepackage[
	bookmarks=true,
	plainpages=false,
	linktocpage,
	colorlinks=true,
	citecolor=green!80!black,
	linkcolor=red!70!black,
	filecolor=magenta,
	urlcolor=magenta,
	breaklinks,
	pdfauthor={Martin Winter},
]{hyperref}

\usepackage{amsmath,amsthm}
\usepackage{calc, mathtools} 

\iftrue
\usepackage{amssymb} 
\else
\usepackage[charter,cal=cmcal]{mathdesign}
\fi

\usepackage{colortbl,color} 
\usepackage[dvipsnames]{xcolor}
\usepackage{centernot} 
\usepackage{array} 
\usepackage{enumitem,moreenum} 
\usepackage{cite} 
\usepackage[nameinlink,capitalize,noabbrev]{cleveref}
\usepackage{nicefrac}

\usepackage{tikz-cd} 

\usepackage[font=small,labelfont=bf]{caption}

\usepackage{blkarray}

\usepackage{inconsolata}
\usepackage{dsfont}

\newcommand{\RR}{\mathbb{R}}    
                   
\newcommand{\Tsymb}{\top}
\newcommand{\T}{^{\Tsymb}}

\def\^#1{^{(#1)}}
\def\s^#1{^{\smash{(#1)}}}

\def\:{\colon}

\definecolor{NiceBlue}{rgb}{0.15, 0.2, 0.75}

\newcommand{\labelstyle}[1]{\upshape(\textit{#1})}
\newcommand{\mylabel}{\labelstyle{\roman*}}

\def\itm#1{{\labelstyle{\romannumeral#1\relax}}}

\newcommand{\nls}{\nolinebreak\space}

\newcommand{\freespace}{\kern.07em}

\newcommand{\bs}[1]{\boldsymbol{#1}}
\newcommand{\bsdot}[1]{\dot{\boldsymbol{#1}}}

\theoremstyle{plain}  
\newtheorem{theorem}{Theorem}

\theoremstyle{definition} 

\crefname{theorem}{Theorem}{Theorems}
\crefname{proposition}{Proposition}{Propositions}
\crefname{lemma}{Lemma}{Lemmas}
\crefname{corollary}{Corollary}{Corollaries}
\crefname{remark}{Remark}{Remarks}
\crefname{example}{Example}{Examples}
\crefname{definition}{Definition}{Definitions}
\crefname{problem}{Problem}{Problems}
\crefname{observation}{Observation}{Observation}
\crefname{construction}{Construction}{Construction}

\DeclareMathOperator{\Span}{span}

\DeclareMathOperator{\vol}{vol}  	
\DeclareMathOperator{\Int}{int}

\DeclareMathOperator{\Div}{div}

\let\eps=\epsilon

\let\<=\langle
\let\>=\rangle
\let\x=\times

\newcommand{\dx}{\,\mathrm dx}

\def\...{...}
\newcommand{\shortStyle}{\textit}
\newcommand{\ie}{\shortStyle{i.e.,}}

\makeatletter
\renewcommand*{\eqref}[1]{%
  \hyperref[{#1}]{\textup{\tagform@{\ref*{#1}}}}%
}
\makeatother

\begin{document}


\expandafter\title
{A short polyhedral proof of the infinitesimal Stoker conjecture}
		
\author[M.\ Winter]{Martin Winter$^\dagger$}
\address{$^\dagger$Max-Planck Institute for Mathematics in the Sciences, Inselstraße 22, 04103 Leipzig, Germany}
\email{martin.winter@mis.mpg.de}
\date{\today}
	
		
\date{\today}
\begin{abstract}
The infinitesimal Stoker conjecture states that an infinitesimal deformation of a convex polytope that preserves all dihedral angles also preserves all normal directions up to isometry.
We give a short polyhedral proof, self-contained up to assuming spectral and geometric properties of the Izmestiev matrix.
\end{abstract}

\maketitle

\section{Introduction}

In 1968 Stoker asked whether the dihedral angles of a convex Euclidean polyhedron determine its face angles \cite{stoker1968geometrical}. 
A first infinitesimal analogue was formulated by Schlenker \cite{m2000dihedral}, which subsequently became known as the \emph{infinitesimal Stoker conjec\-ture}. In this note we work with the following formulation:


\begin{theorem}\label{inf_Stoker}
    Let $P\subset\RR^d,d\ge 2$ be a convex polytope with facets $F_i$, outward unit normal vectors $n_i$ and dihedral angles $\theta_{ij}$.
    Given an infinitesimal deformation of $P$ with $\dot\theta_{ij}\ge 0$ when\-ever $F_i$ and $F_j$ are adjacent, then there exists a skew-symmetric matrix $S\in\RR^{d\x d}$ with $\dot n_i = S n_i$ for all $i$.
    %
\end{theorem}

A first proof of the infinitesimal Stoker conjecture by Mazzeo \& Montcouquiol appeared in 2011 \cite{mazzeo2011infinitesimal} and established the conjecture in dimension three using deformation theory of cone-manifolds.
Another approach in dimension three, though~not~explicitly stated as such, can be extracted from the works of Weiss \cite{weiss2005local,weiss2013deformation}.
Adipra\-sito observed the arbitrary-dimensional infinitesimal statement and sketched an~argument via Hodge–Riemann relations \cite{adiprasito2020geometry}.

In this note we give the first purely polyhedral proof of the infinitesimal Stoker conjecture, and the first direct proof in Euclidean space of arbitrary~dimension.
The proof was initially based on the prestress stability of coned polytope frameworks established in \cite{pachyli2026second}.
It was subsequently distilled to be self-contained~modulo~known spectral and geometric properties of the \textit{Izmestiev matrix} \cite{izmestiev2010colin}.
%
%
The~\mbox{only~ingredients} left hidden within \cite{izmestiev2010colin} are Minkowski's second inequality as well as Bol's condition~for equality. 
For the interested reader, \cref{sec:proof} contains a few comments on how this proof came about, as well as pointers to literature that develops the relevant ideas in greater detail.
%


\section{Setup and assumptions}


All dotted variables like $\dot n_i$ and $\dot\theta_{ij}$ denote formal first-order variations, not~derivatives.
A rigorous treatment can be implemented via dual numbers $\RR[\eps]/\eps^2$. We suppress~such details for readability.
This section lists the implicit assumptions made in the statement of \cref{inf_Stoker}.

    The $n_i$ stay normalized throughout the deformation, hence
    \begin{equation}\label{eq:normalized}
        \<n_i,\dot n_i\> = 0,\quad\text{for all facets $F_i$}.
    \end{equation}
    
    As observed by Karim Adiprasito in \cite{adiprasito2020geometry}, the one-sided condition $\smash{\dot\theta_{ij}\ge 0}$ is easily reduced to its equality version: the classical Schläfli formula states
    $$\sum_{i\ne j} V_{ij}\dot \theta_{ij} = 0,$$
    where $V_{ij}$ is the volume of the ridge $R_{ij}=F_i\cap F_j$.
    Thus, assuming $\dot\theta_{ij}\ge 0$ and by $V_{ij}>0$, we necessarily have $\smash{\dot \theta_{ij}=0}$.
    In the following we will directly assume the equality variant. Since $\theta_{ij} = \pi-\arccos\<n_i,n_j\>$, we can do so by imposing
    \begin{equation}\label{eq:angle}
        \<n_i,\dot n_j\> + \<\dot n_i, n_j\> = 0,\quad\text{whenever $i\sim j$},
    \end{equation}
    where $i\sim j$ means that the facets $F_i$ and $F_j$ are adjacent.
    %
    
    Lastly, the deformations considered preserve the combinatorial type of the polytope. 
    We make use of this by assuming that Minkowski balancing identities on the facets stay valid under variation: for each $i$
    \begin{equation}
        \label{eq:Vn_Vn}
        \sum_{\mathclap{j:j\ne i}} V_{ij}n_{ij} = 0 
        \quad\implies\quad 
        \sum_{\mathclap{j:j\ne i}}\dot V_{ij}n_{ij} + \sum_{\mathclap{j:j\ne i}}V_{ij}\dot n_{ij} = 0,
    \end{equation}
    where $n_{ij}$ is the normal vector of the ridge $R_{ij}$ considered as a facet of $F_i$. We have
    \begin{equation}
        \label{eq:nij}
        n_{ij} 
        := 
        \frac{n_j + \cos(\theta_{ij})n_i}{\sin(\theta_{ij})}
        =
        \frac{n_j - \<n_i,n_j\>n_i}{\sin(\theta_{ij})},
    \end{equation}    
    and since $\dot\theta_{ij}=0$, the variation can be expressed as
    \begin{equation}
        \label{eq:dot_nij}
        \dot n_{ij} 
        = 
        \frac{\dot n_j + \cos(\theta_{ij})\dot n_i}{\sin(\theta_{ij})}
        =
        \frac{\dot n_j - \<n_i,n_j\>\dot n_i}{\sin(\theta_{ij})}.
    \end{equation}   
    


\section{Proof of \cref{inf_Stoker}}


We may assume that $0\in\Int (P)$.
Let $h_i>0$ be the height of the facet $F_i\subset P$ over the origin.
Then the polar dual $Q:=P^\circ$ has vertices $q_i := n_i / h_i$.
We write $\bs q$ for the $(n\x d)$-matrix with rows $q_i$.

The presented proof is in large parts a study of the quadratic form defined by the Izmestiev matrix.%
\footnote{In \cite{pachyli2026second} this corresponds to the energy of the Wachspress stress $\omega^{\mathrm{W}}$. Since we fix the cone~vertex at the origin, it disappeared from all our computations, which simplified the presentation.}
The latter is the only major external ingredient to the proof that we need to assume and that we recall now.

\subsection{The Izmestiev matrix}

Let $M$ be the \emph{Izmestiev matrix} of $Q$. In essence, it is a Colin de Verdière matrix of $Q$'s edge graph. We will use the following properties (see \cite{izmestiev2010colin} or \cite[Theorem 3.3]{pachyli2026second}):
\begin{itemize}
    \setlength{\itemsep}{0.5ex}
    \item[(M0)] $M$ is symmetric.
    \item[(M1)] $M_{ij}=0$ whenever $i\not\sim j$ and $i\ne j$.
    \item[(M2)] $\ker M=\Span\bs q$. 
    \item[(M3)] $M$ has a unique positive eigenvalue.
\end{itemize}
The off-diagonal entries can be given explicitly:
\begin{equation}
    \label{eq:Mij}
    M_{ij} = \frac{V_{ij}}{\|q_i\|\|q_j\| \sin(\theta_{ij})} = \frac{V_{ij} h_ih_j}{\sin(\theta_{ij})}. 
\end{equation}
%
%
Recall that $V_{ij}=\vol_{d-2}(F_i\cap F_j)$. Conveniently for us, $V_{ij}=0$ whenever $i\not\sim j$.

Diagonal entries can be computed as follows: from (M2) we get $M\bs q=0$, or~equi\-valently, $\sum_j M_{ij} q_j=0$ for every $i$. 
Taking the inner product with $n_i$ yields
\begin{equation}
    \label{eq:Mii}
    \sum_j M_{ij} \<n_i,q_j\> = 0
    \quad\implies\quad 
    M_{ii} = -h_i^2 \sum_{\mathclap{j:j\ne i}} \frac{V_{ij}\<n_i,n_j\>}{\sin(\theta_{ij})}.
\end{equation}

We also need the following two facts that follow from straightforward geometric computations (see also \cite[Proposition 3.5]{winter2024rigidity}): 

\noindent
\begin{minipage}{0.5\textwidth}
\begin{equation}
    \label{eq:sum_Mi}
    \sum_j M_{ij} = (d-1)V_i h_i,
\end{equation}
\end{minipage}
\hfill
\begin{minipage}{0.5\textwidth}
\begin{equation}
    \label{eq:sum_Mij}
    \sum_{i,j} M_{ij} = d(d-1) \vol(P),
\end{equation}
\end{minipage}

\noindent
where $V_i:=\vol_{d-1}(F_i)$.

\subsection{Proof scaffold}

Set $\dot q_i := \dot n_i / h_i$ and write $\bsdot q$ for the $(n\x d)$-matrix with rows $\dot q_i$.
We write $(\dot q)_k\in\RR^n$ for the $k$-th column of $\bsdot q$ (\ie\ the vector of $k$-th coordinates of the $\dot q_i$).
The main work for proving \cref{inf_Stoker} is spent on establishing the~following two facts:\footnote{These two steps parallel the proof of prestress stability in \cite[Section 4.1]{pachyli2026second}, though avoid any use of~rigidity terminology.}
$$
(\mathrm{I})\quad 
\bs 1\T\!M (\dot q)_k=0.
%
\qquad\text{and}\qquad 
(\mathrm{II})\quad \sum_k (\dot q)_k\T M (\dot q)_k=0 
$$
where $\bs1=(1,...,1)$ is the all-ones vector.
We can quickly finish the proof from~here.

First, recall
$$
\bs 1\T\! M\bs 1 = \sum_{i,j} M_{ij} \overset{\eqref{eq:sum_Mij}}= d(d-1)\vol(P)> 0.
$$
By (M3) the positive eigenvalue of $M$ is unique, and hence $M$ is negative semidefinite on the $M$-orthogonal complement $\bs 1^{\bot M}:=\{x \mid \bs 1\T\! Mx=0\}$.
In particular,\nls (I) gives $(\dot q)_k\T M (\dot q)_k\le 0$ for all $k$.
With (II) this turns into $(\dot q)_k\T M (\dot q)_k= 0$.
Since $M$ is semidefinite on $\bs 1^{\bot M}$, this implies $(\dot q)_k\in\ker M$.
From (M2) we~get~$(\dot q)_k\in\Span \bs q$,\nls or equivalently, $\dot q_i=A q_i$ for some matrix $A$. 
Substituting~$q_i=n_i/h_i$ and $\dot q_i=\dot n_i/h_i$ yields $\dot n_i=A n_i$.
It remains to show that $A$ is skew-symmetric.

The following argument is standard.
Let $T=\tfrac12(A\T\!+A)$ be the~symmetric part of $A$.
Then
\begin{align*}
    2\<Tn_i, n_j\> &= \<(A + A\T)n_i,n_j\> 
    \\[-1ex]&= \<An_i,n_j\>+\<n_i, An_j\> = \<\dot n_i,n_j\>+\< n_i,\dot n_j\>\overset{\text{\eqref{eq:angle}}}=0
\end{align*}
whenever $j=i$ or $j\sim i$.
So $Tn_i$ is orthogonal to $n_i$ and to all $n_j$ with $j\sim i$. Since $P$ is full-dimensional and convex, the latter contain a basis of $\RR^d$.
Hence $Tn_i=0$ for all $i$.
But also the $n_i$ contain a basis of $\RR^d$. Hence $T=0$, $A$ has no symmetric part, and is therefore skew-symmetric.

\subsection{Proof of (I)}

We show that $\bs1\T\! M (\dot q)_k=0$ for all $k$ simultaneously:
\begin{align*}
\bs1\T\! M \bsdot q = \sum_{i,j} 1_j M_{ij} \bsdot q_i &= \sum_i \Big(\sum_j M_{ij}\Big) \dot q_i 
\\[-1ex]&\overset{\eqref{eq:sum_Mi}}= 
(d-1)\sum_i V_i h_i \dot q_i 
= 
(d-1)\sum_i V_i \dot n_i.
\end{align*}
Consider the tensor field $X_i(x):=x\otimes \dot n_i$. Since $\<n_i,\dot n_i\>=0$, its divergence tangential to $F_i$ is $\Div_{F_i}(X_i)=\dot n_i$.
Applying the divergence theorem in $(*)$ yields\footnote{This argument is a highly compressed application of the vector-valued Schläfli formula from \cite{SS2}.}
\begin{align*}
V_i \dot n_i 
= 
\int_{F_i} \!\!\Div_{F_i}(X_i) \dx
&\overset{(*)}=
\int_{\partial F_i} \!
X_i(x)\cdot \mathrm d n 
= 
\sum_{j:i\sim j} \int_{R_{ij}} \!\!\!\!x\<\dot n_i,n_{ij}\>\dx
\\[-1ex]&= \sum_{j:i\sim j} I_{ij} \<\dot n_i,n_{ij}\>,
\quad\text{where $I_{ij}:=\int_{R_{ij}} \!\!\!\!x\dx$}.
\end{align*}
%
Note that $I_{ij}=I_{ji}$.
The proof finishes by summation over all facets:
\begin{align*}
    \sum_i V_i\dot n_i 
    = 
    \sum_{i\ne j}I_{ij} \<\dot n_i, n_{ij}\>
    &= 
    \sum_{i<j}I_{ij} (\<\dot n_i, n_{ij}\> + \<\dot n_j,n_{ji}\>)
    \\[-2ex]
    \text{\small{subst.\ \eqref{eq:nij} and \eqref{eq:dot_nij}\;$\longrightarrow$}}\quad 
    &= 
    \sum_{i<j}\frac{I_{ij}}{\sin(\theta_{ij})}\overbrace{(\<\dot n_i, n_j\> + \<n_i,\dot n_j\>)}^{=0\text{ by \eqref{eq:angle}}}
    = 0.
\end{align*}

\subsection{Proof of (II)}

The goal is to show $\sum_k (\dot q)_k\T M (\dot q)_k=0$.
For each facet $F_i$,\nls take the inner product of \eqref{eq:Vn_Vn} with $\dot n_i$. Then sum over all facets. This yields
\begin{align}
    \sum_{i\ne j} V_{ij} \<\dot n_{ij},\dot n_i\> 
    &= \notag
    -\sum_{i\ne j} \dot V_{ij} \<n_{ij},\dot n_i\>.
\end{align}
We will rewrite the left side to $\sum_k (\dot q)_k\T M (\dot q)_k$ and show that the right side vanishes.

On the left side we compute
%
\begin{align*}
\sum_{i\ne j} V_{ij} \<\dot n_{ij},\dot n_i\> 
    &\overset{\text{\eqref{eq:dot_nij}}}=
    \sum_{i\ne j} \frac{V_{ij}}{\sin(\theta_{ij})} \big\< \dot n_j - \< n_i,n_j\> \dot n_i,\dot n_i\big\>
    \\&=
    \sum_{i\ne j} \frac{V_{ij}}{\sin(\theta_{ij})} \big(\< \dot n_j,\dot n_i\> - \< n_i,n_j\> \|\dot n_i\|^2\big)
    \\&=
    \sum_{i\neq j} \frac{V_{ij}}{\sin(\theta_{ij})} \<\dot n_i,\dot n_j\> - \sum_i \sum_{j\sim i} \frac{V_{ij}\<n_i,n_j\>}{\sin(\theta_{ij})} \|\dot n_i\|^2
    \\
    \text{\small{subst.\ \eqref{eq:Mij} and \eqref{eq:Mii}\;$\longrightarrow$}}\quad 
    &=
    \sum_{i\ne j} M_{ij}\<\dot q_i,\dot q_j\> + \sum_{i} M_{ii} \|\dot q_i\|^2
    \\&=
    \sum_{i,j} M_{ij}\<\dot q_i,\dot q_j\>
    =\sum_k (\dot q)_k\T M (\dot q)_k.
\end{align*}
And on the right side we compute
\begin{align*}
\sum_{i\ne j} \dot V_{ij} \<n_{ij},\dot n_i\> 
    &\overset{\text{\eqref{eq:nij}}}= 
    \sum_{i\ne j} \frac{\dot V_{ij}}{\sin(\theta_{ij})} \big\<n_j - \<n_i,n_j\> n_i, \dot n_i\big\>
    \\[-2ex]&= 
    \sum_{i\ne j} \frac{\dot V_{ij}}{\sin(\theta_{ij})} \big(\<n_j,\dot n_i\> - \<n_i,n_j\>\overbrace{\<n_i, \dot n_i\>}^{\mathclap{=0\;\text{by \eqref{eq:normalized}}}} \big)
    \\&= 
    \sum_{i\ne j} \frac{\dot V_{ij}}{\sin(\theta_{ij})} \<n_j,\dot n_i\>
    \\[-2ex]&= 
    \sum_{i<j} \frac{\dot V_{ij}}{\sin(\theta_{ij})} \overbrace{\big(\<n_i,\dot n_j\>+\<n_j,\dot n_i\>\big)}^{=0\;\text{by \eqref{eq:angle}}} = 0.
\end{align*}

\iftrue 

\section{How the proof came about}
\label{sec:proof}



In September 2026, correspondence with Arseniy Akopyan revealed that \textit{Wachs\-press Geometry} (in particular, the study of the Izmestiev matrix) allows one to give a polyhedral proof of the (finite) Stoker conjecture \cite{adiprasito2026stokerblog}. 
This new proof builds on techniques developed by the author in \cite{winter2024rigidity}.
Due to independent interest from the side of bar-joint rigidity,
the central ideas of \cite{winter2024rigidity} have been developing towards an infinitesimal setting already since 2023.
This included collaborations with Steven Gortler, Louis Theran and Robert Connelly \cite{connelly2024stress}, as well as Eleni Pachyli and Roman Prosanov \cite{pachyli2026second}.
It was natural to ask whether this line of work might permit~a~likewise elementary approach to the infinitesimal Stoker conjecture.

It quickly became clear that a proof could indeed be extracted from the recent result that coned polytope frameworks are prestress stable \cite[Theorem 1.1]{pachyli2026second}.
%
%
%
The key step was the resolution of the weak stress-flex conjecture, which in turn was based on a vector-valued Schläfli formula first proved by Schlenker \& Souam \cite{SS2}.\nls 
The known proof of the latter did however involve techniques that are not discrete-geometric \cite{SS2}.
Importantly then, one of the contributions of \cite{pachyli2026second} was precisely~to~provide such a discrete-geometric proof for the formula.
At last, a completely polyhedral proof of the infinitesimal Stoker conjecture appeared possible.

Subsequent investigations revealed that the derivation of the Schlenker-Souam-Schläfli formula could be significantly shortened, and for the purpose used in the proof here, even reduced to a simple application of the divergence theorem. This motivated the author to produce the condensed proof presented here.

\fi


\small 

\subsection*{Acknowledgments.}
This short proof is the result of numerous past collaborations, most notably, with Roman Prosanov, Eleni Pachyli, Steven Gortler,  Louis Theran and Robert Connelly.
I am also grateful for their input on earlier drafts of this article.
I am thankful to Arseniy Akopyan, who brought to my attention that my past work contains a proof of the (finite) Stoker conjecture, which triggered this development, and to Karim Adiprasito, who provided exposition through his blog post.

\subsection*{Funding}The author is funded by the SPP 2458 ``Combinatorial Synergies'' (project ID 539851419), funded by the Deutsche Forschungsgemeinschaft (DFG, German Research Foundation).

\subsection*{Statement on AI use}
Generative AI was used during the development and preparation of this article. It assisted in checking calculations, identifying gaps and possible simplifications in preliminary versions of the argument. Some elementary reformulations and shortcuts in the final proof arose from these interactions. 
The overall proof strategy and mathematical development are the author's.
The entire article was written by the author, and all arguments and references appearing in the article were independently checked by the author, who takes full responsibility for the mathematical content.


\bibliographystyle{abbrv}
\bibliography{literature}

\end{document}